\documentclass[12pt]{amsart}
\usepackage{amsfonts,amssymb,amscd,amsthm,amsmath,euscript}
\usepackage{graphpap}

 
\newtheorem{theorem}[subsection]{Theorem}
\newtheorem*{theorem*}{Theorem}

\theoremstyle{definition}
\newtheorem{definition}[subsection]{Definition}

\newtheorem{example}[subsection]{Example}
\theoremstyle{remark}
\newtheorem{remark}[subsection]{Remark}

\theoremstyle{definition}

\newcommand{\mt}[1]{\operatorname{#1}}

\newcommand{\QQ}{{\mathbb Q}}

\newcommand{\OO}{{\mathcal O}}
\newcommand{\RR}{{\mathbb R}}
\newcommand{\PP}{{\mathbb P}}

\newcommand{\NN}{{\mathbb N}}
\newcommand{\FFF}{{\mathbb F}}

\newcommand{\Supp}{\mt{Supp}}
\newcommand{\Pic}{\mt{Pic}}

\newcommand{\NE}{\overline{\mt{NE}}}

\title{$\QQ$-complements on log surfaces}
\author{I.~Yu.~Fedorov and S.~A.~Kudryavtsev}

\date{}

\address{Igor Fedorov: Number Theory Department, Steklov Institute of Mathematics of RAS, 119991 Moscow, Russia}

\email{ifedorov@mi.ras.ru}

\address{Sergey Kudryavtsev: Department of Algebra, Faculty of Mathematics,
Moscow State University, 117234 Moscow, Russia}

\email{kudryav@mech.math.msu.su}

\begin{document}
\thanks{This work was done with the partial support of the
Russian Foundation for Basic Research (grant no. 02-01-00441), the Leading
Scientific Schools (grant no. 00-15-96085) and INTAS-OPEN
(grant no. 2000\#269.)}.

\maketitle
\section*{\bf Introduction}
In this paper  the log surfaces without $\QQ$-complement are classified.
In particular, they are non-rational always. This result takes off the restriction
in the theory of complements and allows one to apply it in the most wide class of log surfaces
$(S,D)$, where the pair $(S,D)$ is log canonical and the divisor $-(K_S+D)$ is nef.
For more information see the papers \cite{Sh} and \cite{Kud}, especially
\cite[Theorems 2.3 and 4.1]{Sh}, \cite[Theorems 2.1 and 3.1]{Kud}.
\par
The work has been completed during the stay at Max-Planck-Institut f$\ddot {\mt u}$r Mathematik in 2003.
We would like to thank MPIM for hospitality and support.

\section{\bf Classification theorem}

We work over an algebraically closed field $k$ of characteristic zero.
The main definitions, notations and notions used in the paper are
given in \cite{Koetal}, \cite{PrLect}.

\begin{definition}
Let $(X,D)$ be a pair, where $D$ is a subboundary. Then a {\it $\QQ$-complement} of
$K_X+D$ is a log divisor $K_X+D'$ such that $D'\ge D$, $K_X+D'$ is log canonical and
$n(K_X+D')\sim 0$ for some $n\in\NN$.
\end{definition}

\begin{example} \cite[Example 1.1]{Sh} \label{ex}
1) Let $\mathcal E$ be an indecomposable vector bundle of rank two and degree 0 over
an elliptic curve $Z$. Then $\mathcal E$ is a nontrivial extension
$$
0\longrightarrow\OO_Z\longrightarrow\mathcal
E\longrightarrow\OO_Z\longrightarrow 0,
$$
see \cite{Hart}. Consider the ruled surface $X=\PP_Z(\mathcal E)$. Let $C$ be the
unique section corresponding to the exact sequence. Note that $C|_C\sim\det \mathcal E\sim 0$,
$-K_X\sim 2C$ \cite{Hart}. In particular, $-K_X$ is nef and Mori cone $\NE(X)$ is generated by two
rays $R_1=\RR_+[C]$ and $R_2=\RR_+[f]$, where $f$ is a fiber.
\par
We claim that $C$ is the unique curve in $R_1$. In fact, let there is a curve
$L\ne C$ in $R_1$ then $L\equiv mC$, where $m=L\cdot f\ge 2$. It is easy to prove that
$L\sim mC$.
Therefore the linear system $|L|$ gives a structure of an elliptic fibration on $X$ with a
multiple fiber $C$. Hence $C|_C$ is an $m$-torsion element in $\Pic(C)$, a contradiction
with $C|_C\sim 0$. We proved that the log surface $(X,cC)$ does not have $\QQ$-complement,
where $0\le c\le 1$.

\par
2) Consider the pair $(X,C)$ from the previous example.
Let $P_1,\ldots,P_r$ be the arbitrary points of $C$. Take any number of
blow-ups at $P_1,\ldots,P_r$. We obtain the pair $(\widetilde X, \widetilde C)$, where
$K_{\widetilde X}+\widetilde C=g^*(K_X+C)$ and $\widetilde C$ is a proper transform of
$C$. It is clear that  $(\widetilde X, \widetilde C)$ does not have $\QQ$-complement.
If we contract any chains of $(-2)$-curves on $\widetilde X$ and maybe $\widetilde C$ (if it is possible)
then the log surface
obtained does not have $\QQ$-complement also.

It is obvious that all these log surfaces do not have complement too
(see the definition of complement in \cite[Definition 4.1.3]{PrLect}).
\end{example}

\begin{theorem}
Let $S$ be a normal projective surface and $D$ be a boundary on $S$ such that
$K_S+D$ is log canonical and $-(K_S+D)$ is nef. Assume that $(S,D)$ does not
have $\QQ$-complement. Then the pair $(S,D)$ is of example $\ref{ex}$, in particular
$S$ is non-rational.
\begin{proof}
Let $f\colon \widetilde S\to S$ be a minimal resolution and
$K_{\widetilde S}+\widetilde D=f^*(K_S+D)$. Then the pair $(\widetilde S, \widetilde D)$
does not have $\QQ$-complement too. By abundance theorem \cite[Theorem 8.5]{Koetal}
the kodaira dimension $k(\widetilde S)=-\infty$. Consider two cases.
\par
A) Let $\widetilde S$ be a rational surface. Since $\widetilde S\not\cong \PP^2, \FFF_0$
then some model of $\widetilde S$ is $\FFF_n\  (n\ge 1)$. We have
$g\colon \widetilde S\to\FFF_n\to Z\cong \PP^1$. Let $\widetilde E_{\infty}$ be the
proper transform of the minimal section of $\FFF_n$.
\par
Now we construct the divisor $\widetilde D'\ge \widetilde D$ such that
$K_{\widetilde S}+\widetilde D'\equiv 0/Z$,
$(K_{\widetilde S}+\widetilde D')\cdot \widetilde E_{\infty}=0$ and the pair
$(\widetilde S, \widetilde D')$ is log canonical.
Hence $K_{\widetilde S}+\widetilde D'\equiv 0$ and $K_{\widetilde S}+\widetilde D'$ is
a $\QQ$-complement of $K_{\widetilde S}+\widetilde D$ by
abundance theorem, a contradiction.
\par
Let $f_1, \ldots, f_k$ be the reducible fibers of $g$. Let $f'_1, \ldots, f'_k$ be any
irreducible components of $f_1, \ldots, f_k$. By considering the linear system
$|mE_0|$ on $\FFF_n$, where $E_0$ is a zero section and $m\gg 0$ we obtain a free pencil
$|L|$ on $\widetilde S$ such that $L\cdot \widetilde E_{\infty}=0$ and
$f_i\cap L_{\mt{gen}}\subset f'_i$ are the generic points of $f'_i$ for the general
element $L_{\mt{gen}}\in |L|$. Hence adding the required number of different
$L_{\mt{gen}}$ and general fibers of $g$ we obtain $\widetilde D'$.
\par
B) Let $\widetilde S$ be a non-rational surface. Then we have a contraction onto a curve
$f\colon \widetilde S\to Z$, where general fiber is $\PP^1$ and $p_a(Z)\ge 1$. By
\cite[Lemma 8.2.2, Corollary 8.2.3]{PrLect} no components of $\Supp \widetilde D$
are contained in the fibers of
$f$ and the pair $(\widetilde S, \widetilde D)$ is canonical. There are two variants.
\par
1) Let $\widetilde D=\widetilde C+\widetilde D_1$, where
$\widetilde C$ is an irreducible curve. Hence $p_a(Z)=1$ and we have a birational contraction
$\widetilde S\to \overline S$ such that $K_{\widetilde S}+\widetilde D\equiv 0/\overline S$,
where $\overline S$ is a minimal model of $\widetilde S$. Therefore the pair
$(\overline S, \overline D)$ is without $\QQ$-complement, where $\overline D$ is the image
of $\widetilde D$. By \cite[Corollary 2.2]{Sh} and \cite{Hart} the surface $\overline S$ is
of example \ref{ex} 1) and $\overline D=\overline C$. Q.E.D.
\par
2) Let $\widetilde D=\sum d_i\widetilde D_i$, where $d_i<1$ and let $\overline S$ be a minimal
model of $\widetilde S$. By \cite[Proof of theorem 3.1]{Kud} we have $E^2\ge 0$ for every curve
$E$ on $\overline S$, $p_a(Z)=p_a(\overline D_i)=1$,
$K_{\overline S}\cdot \overline D_i=\overline D_i^2=\overline D_i\cdot\overline D_j=0$
for all $i,j$.
It follows easily that we can consider $\widetilde D'=\widetilde D_1+\sum_{i\ge 2}d_i\widetilde D_i$
instead of $\widetilde D$ by \cite{Hart} (if $\widetilde D=0$ then $\widetilde D'=\widetilde C$,
where $\widetilde C$ is a proper transform of an irreducible curve $\overline C$ such that
$K_{\overline S}\cdot\overline C=\overline C^2=0$).
The problem is reduced to the variant 1).
\end{proof}
\end{theorem}

\begin{remark}
The log surface without complement is the same one. This fact is proved similarly.
Thus the log surfaces with complement are equivalent to the log surfaces with
$\QQ$-complement.
\end{remark}

\end{document}